# Estimating impacts of covid19 on transport capacity in railway networks

**Nikola Bešinović[1,*], Christopher Szymula[1,2]**

[1]Digital Rail Traffic Lab, Delft University of Technology, The Netherlands

[2]"Friedrich List" Faculty of Transportation and Traffic Science, Technische Universität Dresden, Germany

*Corresponding author: n.besinovic@tudelft.nl

## Abstract
Due to the covid19 crisis, public transport (PT) systems are facing new challenges. Regarding restrictive measures such as physical distancing and the successive returning of passengers after the "intelligent lockdown", significant lack of transport capacity can be expected. In this paper, the transport capacity of a PT network is assessed, using a mathematical passenger route choice and train scheduling model. By analysing the overall number of transported passengers and the resulting link and train utilization; the networks capabilities of facilitating different demands under capacity restrictions (e.g. physical distancing) are addressed. The analysis is performed on the Dutch railway network. The results show that at most 50% of the pre-covid19 demand can be transported, while most of the trains will be highly utilized reaching their maximum occupation. Thus, significantly more parts of the network are becoming highly utilized, leading to a more congested and vulnerable system than in normal conditions before covid19.

## 1. Introduction

Researchers in public transport (PT) are usually dealing with problems like overcrowded trams in peak hours, delayed trains, or other side effects of mobility in highly utilized transportation systems. But due to the covid19 crisis, we as mobility professionals are no longer concerned with these unintended side effects of mobility. In fact, public transport now only transports 10% of the usual passengers and there is relatively more cycling and walking than ever.  For example, some of the largest decrease in use of public transport has been seen in San Francisco Bay Area - 95%[1], Lyon - 92%[2], Rome - 90%[3]. In the Netherlands, decrease in use in peak hours is about 90-95% as well[4]. After pandemic passes and the "intelligent lockdown" is lifted, some restrictive measures are expected to stay much longer. In particular, we already observe that a physical distancing of 1.5m is becoming our new reality. We shall then accept that with returning of passengers to the PT system, transport capacity (number of seats) will not be sufficient in peak hours. Therefore, infrastructure managers and operators need to reorganize their services in response to these new societal conditions.

---

[1] https://transitapp.com/coronavirus, visited 25 May 2020.
[2] https://transitapp.com/coronavirus, visited 25 May 2020.
[3] https://moovitapp.com/insights/en/Moovit_Insights_Public_Transit_Index-countries, visited 25 May 2020.
[4] https://www.dat.nl/nvp/, visited 25 May 2020.

New challenging questions arise: **What does peak traffic look like in public transport in such a society? What is the transport capacity of our system under physical distancing?** In particular, we shall understand what the system capabilities are today, whether we have available capacity to transport all passenger demand and if not, how much demand we can still satisfy with the current railway system.

Bešinović (2020) defined resilience of railway transport systems as the ability to withstand as well as recover quickly from various disruptions such as internal resource failures to weather-related consequences. However, the current epidemics could be added to the list of external threats as it represents a new type of disruption that has not been addressed in the literature previously. What is more, impacts of epidemics are expected to be visible longer such as physical distancing and time needed for recovering demand to levels of pre-covid19 period (van Ammelrooy, 2020). Therefore, our railway system needs to become resilient against these new challenges as well. To do so, the first steps would be to assess the new system capacity due to covid19. However, immediate negative effects on reduced system performance due to limited resources caused by covid19 (e.g. reduced number of seats), has not been considered so far. In this paper, we address railway transport capacity **during capacity-related disruptions being** extreme shortage as a consequence of covid19 and focus on estimating the maximum number of passengers that can travel during peak hours under various passenger transport demand.

We distinguish between traffic capacity and transport capacity. **Traffic capacity** is defined as the number of trains that can be scheduled within a given period of time on a given infrastructure**.** This could be referred as to infrastructure capacity as well. **Transport capacity** is defined by number of seats offered to passengers within a given period of time. In the existing literature, models for assessing railway traffic capacity have been mostly addressed (Bešinović and Goverde, 2018) and to this purpose various applications have been considered including analytical models, simulations, and optimization models. However, only limited research exists on estimating transport capacity of railway systems.

**In this research, we use a mathematical passenger route choice and train scheduling model to explore transport capacity impacted by covid19 in the Dutch railway network.** The mixed integer programming (MIP) model incorporates operator's and passengers' perspectives, it maximizes number of transported passengers on one side and finds attractive passenger routes through the network on the other side. The model has been efficiently solved using a heuristic based on column generation with mixed integer linear programming to efficiently generate alternative passenger flows. We apply the model on the real-world instance of the Dutch railway network and assess network capacity on link level and train level. The results show that at most 50% of the pre-covid19 demand can be transported, while most of the trains will be highly utilized reaching their maximum occupation.

The main contribution of the paper is threefold.

- We provide an in-depth empirical sensitivity analysis of railway capacity under passenger demand and train capacity variations considering multiple system indicators such as passenger volume, link utilization and train utilization.
- We explore impacts of covid19 on transport capacity of the real-life Dutch railway network.
- We give directions for future research on tackling challenges caused by covid19 and extreme shortage of transport capacity in general.

The remainder of the paper is as follows. Section 2 reviews methods for capacity assessment and covid19-related research. Section 3 discusses new restrictions on vehicle capacity. Section 4 describes the mathematical model for estimating transport capacity of railway networks. Section 5 presents the results of the evaluation of using railway transport under different passenger transport demands. Section 6 gives concluding remarks and points out directions for future research.

## 2. Literature review

In this section, we review literature on traffic and transport capacity assessment and recent research related to covid19 within transport domain. Approaches for capacity assessment typically include analytical models, queuing, simulation and optimization. Some of the most common indicators used for assessing capacity are capacity occupation and capacity utilization, which could be determined for vehicles, stations, lines, corridors, or networks (Vuchic, 2007). Occupation represents the amount of trains in the system or number of occupied seats in the vehicle and represents an absolute number of units. Instead, utilization represents a ratio between used/occupied capacity and the scheduled/offered capacity and ranges typically between 0 and 1 (in some congested systems can result in values greater than 1 as well). By scaling the occupation on the offered capacity, the latter is more suitable for comparing vehicles with different capacities, so we use it in our experiments.

Various approaches for traffic capacity assessment can be found. The UIC 406 capacity method, i.e. compression method, and its extensions analytically analyze the consumed capacity of a given railway infrastructure (UIC, 2013; Landex, 2008). Schwanhäußer (1978, 1994) introduced a queuing theory approach for evaluating the capacity while several extensions of this approach can be found (Büker and Seybold, 2012; Huisman, et al, 2002; Wendler 2007; Weik, et al. 2016). Krueger (1999) and Lai and Barkan (2009) proposed parametric modelling. Optimization models for capacity assessment are presented in e.g. Burdett and Kozan (2006) and Burdett (2015). Quaglietta et al. (2014), and Jensen et al. (2017) deployed simulations for this purpose. Based on extensive experience, it has been concluded that timetable structures are required to understand the complex interactions in a dense and complex railway network (Odijk et al., 2006), and Eliasson and Börjesson (2014). Railway traffic capacity depends on various aspects that can be traditionally categorized in three groups: infrastructure, rolling stock, and traffic management. Goverde et al. (2013) showed the influence of various signalling systems on the capacity occupation. Strategic Rail Authority (2014) stressed the importance of traffic management and operational rules like dominant train type (passenger, freight or mixed), use of tracks (unidirectional/bidirectional), mix of train services with different characteristics (speed, stopping pattern, frequency), train sequences, dwell times and connections in stations. UIC (2013) explained that capacity depends on the way the infrastructure is utilized which is represented in the capacity balance of the number of trains, the average speed, the traffic heterogeneity, and stability. A detailed analysis of different aspects affecting capacity can be found in Abril et al. (2008); Harrod (2009); Schmidt (2014); Shih et al. (2014) and Lindfeldt (2015), while an empirical comparison of different capacity assessment methods can be found in Rotoli et al. (2016). For more detailed overview of methods for railway traffic capacity we refer to Bešinović and Goverde (2018).

Gentile et al. (2016) reviewed transit assignment models and addressed both frequency-based and schedule-based variants, and possible approaches to represent the dynamic system such as steady state, macroscopic flows, and agent-based (route choice, flow propagation, arc performances). Cats and Hartl (2016) compared a schedule-based model and an agent-based transit assignment model focusing on the congestion conditions in transit networks during everyday operations. To do so, they used simulation models VISUM and BusMezzo. At a more abstract level, Wu et al. (2015) established an estimation model of urban transportation supply-demand ratio (TSDR) to quantitatively describe the conditions of an urban transport system and to support transport policy-making based on the system dynamic modelling. Using automated fare collection and train tracking data, Hänseler et al. (2020) described the realized pedestrian movements in train stations and vehicle-specific train ridership distributions for a line in the Dutch railway network. Szymula and Bešinović (2020) jointly determined critical infrastructure disruptions and passenger demand that can be transported in order to assess vulnerability of railway networks. The solution approach includes column and row generation to efficiently incorporate passenger paths and timetable adjustments, respectively.

Current studies on traffic and transport modelling commonly considered normal transport conditions, and occasionally changes/failures in physical system resources. While doing so, the focus is most often on infrastructure capacity, uncongested transit assignment and user equilibrium in congested networks, and also on passenger routing and comfort, assuming that all passengers can travel, rather than on estimating transport capacity of the system. However, covid19/pandemics represents a disruption of a different nature and significantly impacts passengers' interactions by requiring a 1.5m physical distancing. As such it poses a much greater restriction to the system performance leading to many cancelled travelling requests and such system assessments have not been tackled so far. To this end, the focus of our paper is on estimating maximum number of passengers that can travel during peak hours. Regarding PT and railway networks, the research on (remaining) available capacity due to covid19 is still missing.

Some researchers already tackled epidemiological spreading of viruses by mathematical modelling of the dynamics (Prasse et al. 2020; Van Mieghem, 2020; Khrapov and Loginova, 2020). Still, only limited research in PT and railway networks can be found. Krishnakumari and Cats (2020) addressed virus spreading using epidemiological modelling and contact networks to evaluate effects in the metro network in Washington DC. Kucharski and Cats (2020) examined how ride-sharing services contribute to exacerbating the epidemic using an epidemiological model. They assume that initially infected travellers infect others before they quarantine, and that demand may vary from one day to the other. There have been a several surveys investigating impacts of covid19 on passenger transport. A recent survey on risk perception by Gerhold (2020) identified that "people especially worry about being infected in places with high public traffic such as public transport and shops or restaurants". In the Netherlands, a survey showed that people generally expect to use public transport less than before covid19 (Jacobs, 2020). However, not all implications are negative to society. Badii et al, (2020) and Sommer (2020) reported that there have been remarkable improvements in air quality and declines in petroleum use throughout the world as a result of the short-term and significant decline in freight and passenger travel. The first steps have been towards designing railway timetables for excessive passenger demand while satisfying the demand as much as possible for given infrastructure (Bešinović et al., 2019). Authors defined the demand implicitly representing it by predefined desired line plan and train line frequencies. The objective was then to schedule as many lines and with high frequencies as possible. Gkiotsalitis and Cats (2020) proposed a mathematical model for determining frequencies of services in metro networks considering physical distancing due to covid19.

## 3. Vehicle capacity

Due to physical distancing measures of 1.5m, passengers would need to keep bigger distances, which essentially means that not all seats in the vehicles can be occupied. Table 1 shows the original and reduced capacity due to covid19. In the Netherlands, currently the train capacity is organized so as one person is allowed to sit per row (i.e. 4 seats). Based on these values we expect that vehicles would be used approximately up to 1/4 of their original seating capacity. In urban systems we see some greater reductions in vehicle capacity. For example, Krishnakumari and Cats (2020) reported the reduced capacity to be around 20% of the original capacity. Note that standing capacity goes from 4-6 pax/m$^2$ to approx. 0.5 pax/m$^2$, which reduces the available capacity even more severely. Therefore, the higher reductions on vehicle capacity may indeed be expected in urban transit systems with a bigger share of standing capacity (metro, trams, LRT) than in conventional railway systems assuming the physical distancing of 1.5m.

*Table 1. Train capacity in some PT systems in normal and covid19 conditions*

| Vehicle type, operator (country) | Original capacity | covid19 capacity |
|---|---|---|
| intercity train, NS (NL) | 400-1100 seats (max capacity for VIRM4+VIRM6) | 100-250 |
| local train, NS (NL) | 320-540 | 80-140 |
| Bus articulated, RET (Rotterdam, NL) | 120 | 30 |
| Bus regular, RET (Rotterdam, NL) | 80-100 | 20 |
| Bus regular, TriMet[5] (Portland, USA | 60-80 | 10-15 (10 individuals or up to 15 if people are riding together (such as couples or parents with children) |
| Metro (6-8 cars), Washington DC, USA (Krishnakumari and Cats, 2020) | 1700 | 312 |

---

[5] https://blog.trimet.org/2020/04/02/covid-19-update-limiting-how-many-people-are-allowed-on-transit/

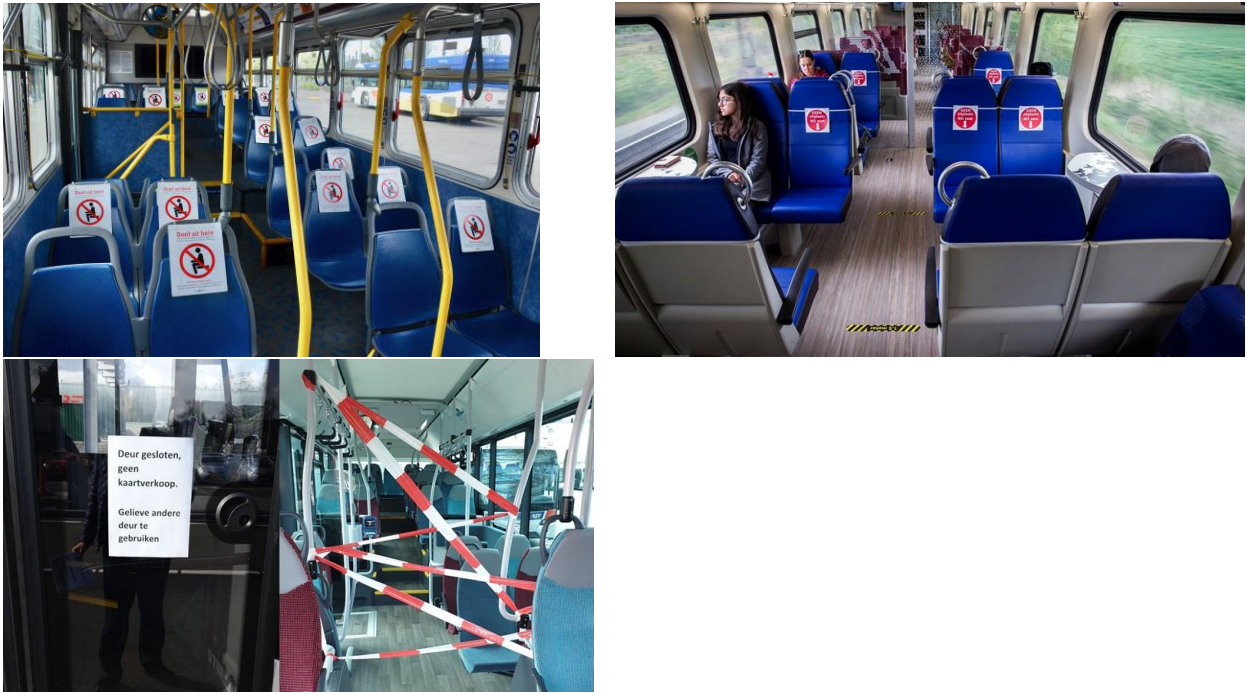

*Figure 1. Seats marked during covid19: TriMet bus[6] (top left), NS intercity train[7] (top right), RET bus[8] (bottom)*

## 4. Modelling

To evaluate transport capacity in railway networks under physical distancing we use an adjusted version of the model by Szymula and Bešinović (2020). In our current model for transport capacity assessment, compared to Szymula and Bešinović (2020), no infrastructure disruptions are considered, trains are considered to run according to schedule without rerouting possibilities, and a simplified objective function is defined. Additionally, input parameters such as train capacity are varied, representing disrupted transport capacity. The model provides the passenger flows for given timetables, train routes and train capacities in a railway network. In particular, the model maximizes the number of transported passengers and routes passengers over shortest available paths. It combines arc-based and path-based formulations to efficiently model trains and passengers, respectively. Therefore, the capacity assessment model addresses both operator's as well as passengers' perspective.

The railway network is modelled as a graph of nodes and links with stations being assigned as nodes $i \in N$ and railway sections as the links $(i,j) \in A$. Trains are denoted by $t \in T$, and train routes are represented by binary variables $x_{ij}^t$ for each train $t$ over link $(i,j)$. To model passenger behaviour, the passenger flows in the network are estimated based on the demand volumes $d_k$ of all Origin-Destination (OD) pairs $k \in K$ and their corresponding OD distances $c_p^k$. By maximizing the passengers traveling on the shortest available paths, the actual transported demand is generated, and the passenger flows are routed correspondingly. While doing so, capacity constraints such as the given train capacities (i.e. number of seats $s^t$) in the traversed connections are regarded. Finally, share of passengers of OD pair $k$ using a path $p$ is denoted as $f_p^k$. The mathematical formulation of the model used for transport capacity assessment is as follows:

$$\max \sum_{k \in K} \sum_{p \in P^k} \frac{d_k}{c_p^k} f_p^k \tag{1}$$

---

[6] By Steve Morgan - Own work, CC BY-SA 4.0, https://commons.wikimedia.org/w/index.php?curid=89201193
[7] https://www.rtvdrenthe.nl/nieuws/159502/NS-test-aangepaste-treinen-op-traject-door-Drenthe
[8] https://www.rijnmond.nl/nieuws/193032/Achter-in-de-RET-bus-stappen-svp

Such that:

$$\sum_{j \in N} x_{ij}^t - \sum_{j \in N} x_{ij}^t = \begin{cases} -o_i^t, & \text{if node } i \text{ is a starting node} \\ d_i^t, & \text{if node } i \text{ is an end node} \\ 0, & \text{else} \end{cases} \quad \forall t \in T, i \in N^t \tag{2}$$

$$\sum_{i \in N^t} o_i^t = \sum_{i \in N^t} d_i^t \qquad \forall t \in T \tag{3}$$

$$o_i^t = 1 \qquad \forall t \in T, i = N(O_t) \tag{4}$$

$$d_i^t = 1 \qquad \forall t \in T, i = N(D_t) \tag{5}$$

$$o_i^t = 0 \qquad \forall t \in T, i \neq N(O_t), N(D_t) \tag{6}$$

$$d_i^t = 0 \qquad \forall t \in T, i \neq N(O_t), N(D_t) \tag{7}$$

$$\sum_{t \in T} x_{ij}^t \leq CAP_{ij} \qquad \forall (i,j) \in A^t \tag{8}$$

$$\sum_{p \in P^k} f_p^k \leq 1 \qquad \forall k \in K \tag{9}$$

$$\sum_{k \in K} \sum_{p \in P^k} \delta_{i,j}^p d_k f_p^k \leq \sum_{t \in T} s^t x_{ij}^t \qquad \forall (i,j) \in A \tag{10}$$

$$x_{ij}^t, o_i^t, d_i^t \in \{0,1\} \qquad \forall t \in T, n \in N^t \tag{11}$$

$$f_p^k = [0,1] \qquad \forall k \in K, p \in P^k \tag{12}$$

Equation (1) describes the objective of the optimization problem, which is to maximize the shares $f_p^k$ of the passenger demand $d_k$ for all OD-pairs $k \in K$ on paths $p \in P^k$ at minimal cost. The minimization of the path cost is considered by utilizing the inverse path cost $\frac{1}{c_p^k}$ over all existing paths. Equation (8) ensures the flow continuity of the offered train services $x_{ij}^t$ for all trains $t \in T$ in all nodes $i \in N^t$ of the corresponding train route and that trains are only allowed to start (end) at their corresponding starting (terminating) nodes. Equation (9) is balancing the number of source and sink nodes along each train route. Equations (10) and (11) fix the source and sink nodes. Equations (12) and (13) are preventing train services of being terminated outside of source or sink nodes, respectively. Equation (8) sets the maximum infrastructure capacity $CAP_{ij}$ for the train services on each link $(i,j) \in A$ in the network. Equation (9) restricts the passenger flow shares for each OD pair $k$ to be at most as large as the overall demand of that pair. Equation (10) restricts the cumulative passenger flow shares on each arc to be at most as large as the offered accumulated seating capacity $s^t$ of the corresponding train services $x_{ij}^t$ that traverse link $(i,j) \in A$. Parameter $\delta_{i,j}^p$ identifies whether the actual link is traversed by the corresponding passenger-flow path $p$ ($\delta_{i,j}^p = 1$) or not ($\delta_{i,j}^p = 0$) . Equations (11) and (12) are restricting the range of the decision variables.

For tackling large-scale challenges of real-life railway networks, a solution framework has been introduced. In theory, passengers typically have extremely large numbers of travelling path alternatives. And so, considering them all in the model at once may render the problem impossible to solve (e.g. Boulliet et al., 2007 p.128; Gentile et al., 2016, p.296). In practice, many of those alternative paths may not be actually used due to their unattractiveness, e.g. represent excessively long detours. In order to keep our model computationally solvable, we apply a column generation inspired heuristic to identify new beneficial available passenger paths that can be added iteratively to the optimization model (Lusby et al. 2012). Thus, these so called available shortest paths are iteratively added in the column generation approach for dealing with overcrowded trains, satisfiable demand and acceptable detours. The detailed solution framework is explained in Szymula and Bešinović (2020).

## 5. Experimental results
### 5.1 Setup

The Dutch passenger railway network, one of the most utilized railway networks in Europe, has been chosen for performing our experiments to assess transport capability of the railway system. The train-related data is received from the General Transit Feed Specification (GTFS) data, which consists of the train lines, routes and the scheduled arrival and departure times. The passenger demand used is based on realistic demand data from the Dutch railway network. In particular, OD pairs with a demand larger than 100 passengers/hour are considered in our experiments to focus on more significant passenger flows. We assume the original pre-covid19 timetable (with all services running) and put the focus of the sensitivity analysis on evaluating the system performance on changing the passenger demand levels under original and reduced train capacities. In total, 507 trains are considered within peak hour.

We use 5 variants of demand size. First, the original transport demand of normal conditions is used and referred to as 100%. Due to "intelligent lockdown", passenger demand reduces significantly, so we also consider four additional restricted variants with 5%, 25%, 50% and 75% of the original demand. These demand shares can be understood as a gradual increase of demand after the restrictions are lifted. The same demand shares are also considered under covid19 conditions, which cause a significant reduction in provided service capacities due to physical distancing. In the scenarios, the relative demand is uniformly altered for all OD pairs. Such demand changes seem reasonable knowing that the lockdown measures are typically made on a national level and imposed/released nationwide. However, differences may still occur due to different acceptance attitude towards lockdown measures between certain regions and changes in the overall travel behaviour of different passenger groups. In the sensitivity analysis, we evaluate the transported demand, link utilization and train utilization and report it for both normal and covid19 conditions in Sections 5.2-5.4.

For covid19 conditions, based on observations in Section 3 on reduced train capacity in the Dutch railways, we adopt the capacity of trains to be 200 passengers/train. Also, it is assumed that all passengers travel alone. This may be reasonable during peak hours as most passengers commute to their work. In particular, for some intercity lines this could be a minor underestimate, while for some local trains it would be an overestimate. Given that the Dutch timetable includes around approximately the same number of intercity and local trains, this tends to be a good estimate. In addition, capacity of trains in normal conditions is 1000 passengers/train. Finally, in order to simplify routing, the travel cost is calculated per traversed arc using averaged travel cost over all passenger trains running on this arc. This allows to connect passenger paths to served relations rather than single trains.

### 5.2 Transported demand

Figure 2 shows expected transported passengers for a given passenger demand during a peak hour in the entire network. It can be seen, that for very small demands – both covid19 and normal capacity – the network can transport an (almost) equal number of passengers. However, by increasing the demand, the number of transported passengers that cannot be transported due to the reduced capacity in covid19 scenario increases significantly. Once the demand returns to its normal size, the railway system under physical distancing would be capable to transport as many as 165,877 passengers in the peak period. **It can be observed that with 50% of passenger demand share the number of unserved passengers reaches approximately 30% (see 1 in Figure 1) and this gap grows further up to 50% difference (see 2) compared with the original (pre-covid19) demand.** Figure 3 shows passenger flows in the network in normal and covid19 conditions with 50% demand share. It becomes clear that in the latter significantly fewer passengers are being transported.

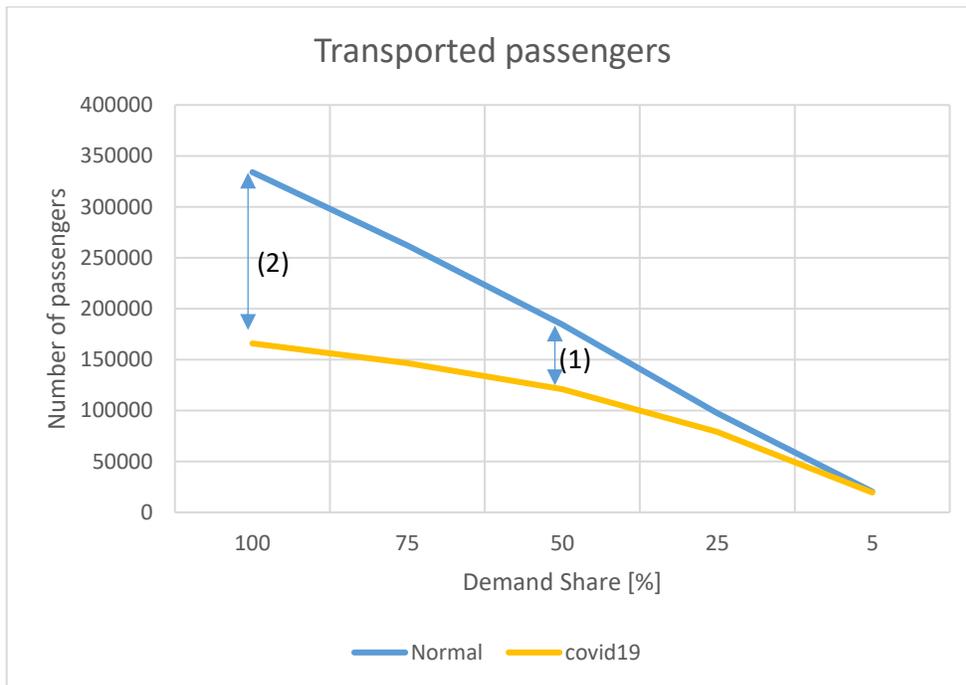

*Figure 2. Transported passengers based on the actual given demand*

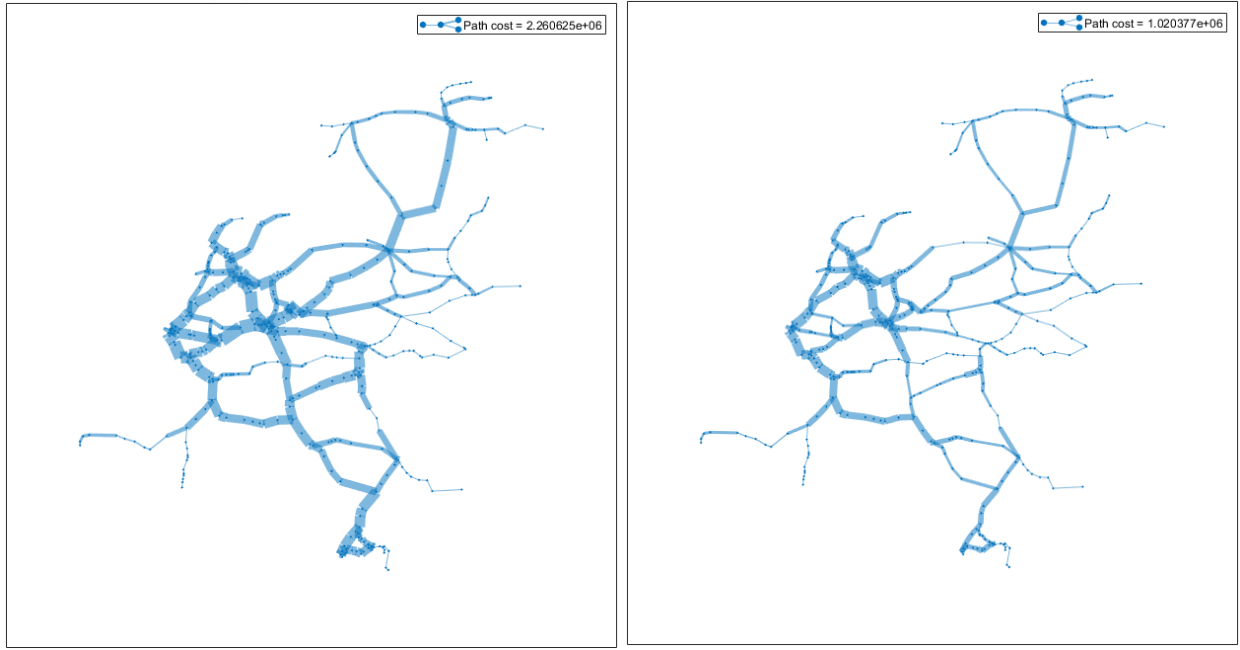

*Figure 3. Passenger flows for normal and covid19 conditions with 50% demand.*

### 5.3 Link Utilization

Figure 4 shows statistics of link utilization under normal and covid19 conditions reporting mean, median and standard deviation (stdDev) for given demand shares. As it can be seen, the average link utilization under normal conditions decreases with lower demands, it ranges from 0.33 to only 0.03. Since the overall number of passengers decreases with lower demand, the overall mean and median utilization decrease continuously as well as the stdDev. It is noticeable, that the median is constantly lower than the mean. This indicates that larger share of links seems to be less frequently utilized. A correspondingly decreasing stdDev, from 0.30 to 0.17, indicates that, as fewer passengers are travelling in the network, the heterogeneity in link utilization decreases as well. Thus, the links are getting utilized more equally. Since the demand is decreasing, a more homogeneous (i.e. lower) utilization likely means transporting a (very) limited number of passengers and thus a lesser utilized network.

It can be also seen that in the covid19 case, a continuously decreasing link utilization comes along with decreasing demand shares, it ranges from 0.49 to 0.12. In addition, for the original demand, the median is larger than the mean utilization, 52% (0.52) compared to 0.49, and drops below the mean only for lower demand shares. This indicates an over-average utilization of links, which could result from the high utilization of major network parts, indicating a potential overall network saturation. As the demand in the network decreases, the link utilization decreases with fewer links being highly saturated, which leads to a bigger decrease of median link utilization. Overall, the network remains with rather high link utilization for all demand scenarios. Like the normal case, the decreasing stdDev illustrates the slow but steady reduction in traveling passengers, but with less significant variations between different demand shares.

Comparing the normal and the covid19 case from a network perspective, they show the railway network performance (link utilization) at different performance levels. Figure 5 shows link utilization under covid19 and normal conditions of 100% demand share. Due to the reduced train capacities in the covid19 case, the network reaches a highly utilized level, showing high link utilizations in major parts of the network. Hence, decreasing demands lead to fewer passengers travelling in the network, but still causing many highly utilized links. The normal case results in a significantly lesser utilized network, where only few areas in the network are highly utilized, while most of the links are used by a limited number of passengers only, and thus leaving significant capacity unused. **In addition, for the covid19 condition and original demand, as much as 20% of all links would be fully utilized, while for the lightest demand share, it is expected to be around 5% of all links**. At the same time, for normal case and original demand, this value would reach maximally about 5% of links.

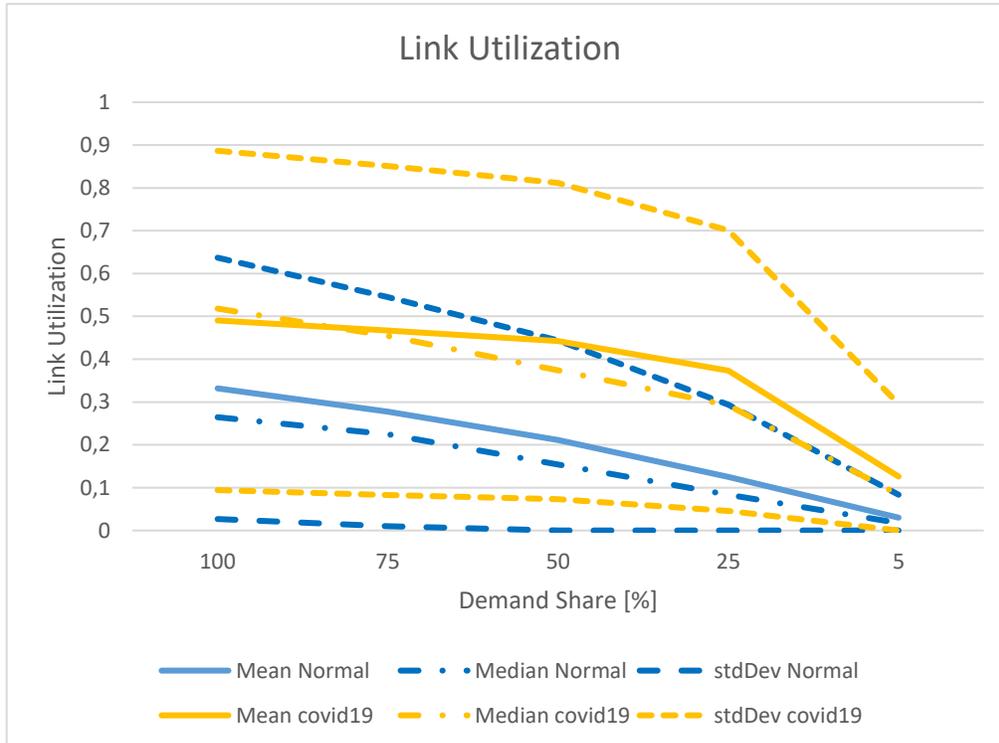

*Figure 4: Link Utilization normal and covid19 case statistical indicators*

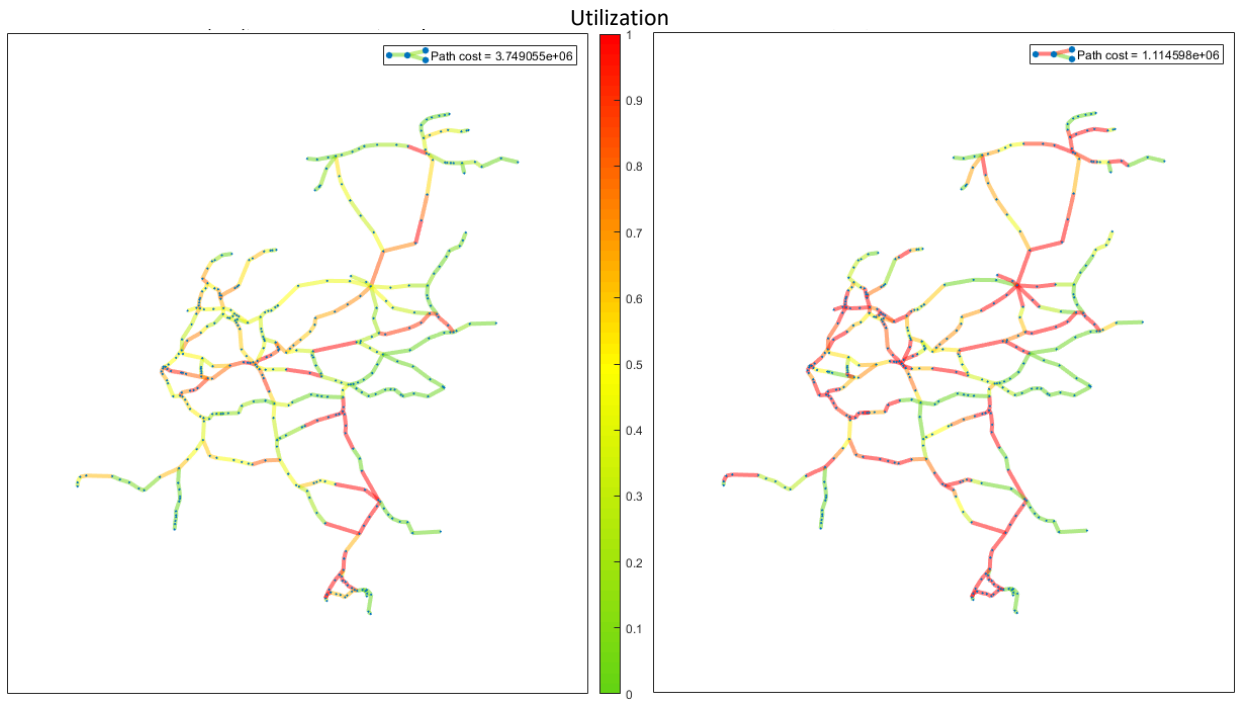

*Figure 5: Link Utilization of 100% Demand Share in normal and covid19 case*

## 5.4 Train utilization

We report statistics on **average train utilization** (average number of passengers onboard on each train over the complete train journey) as well as **maximum train utilization** (number of passengers onboard on each train at the busiest link of the train journey).

Figure 6 shows the average train utilization, its mean and stdDev in the network. During normal conditions, the average train utilization stays remarkably low. For demands smaller than 50%, it is not bigger than 20% and looking at the median, half of all trains have average train utilization of 0.2 and lower. This extends with the fact that 80% of trains transporting at most 160 passengers, while it fluctuates between 100 and 360 passengers.

It becomes clear that due to covid19, an average train utilization would be at 0.40 (i.e. 80 out of 200 seats occupied) for 25% demand share. More interestingly the mean values are always higher than the average, indicating a high utilization of vehicles. With 50% demand share, trains would carry on average 100 passengers, which represents average utilization of around 0.50 of trains capacity. It is also observed that **at least 25% of trains have an average capacity utilization greater than 0.70 on board for scenarios above 50% of normal demand** (Figure 7).

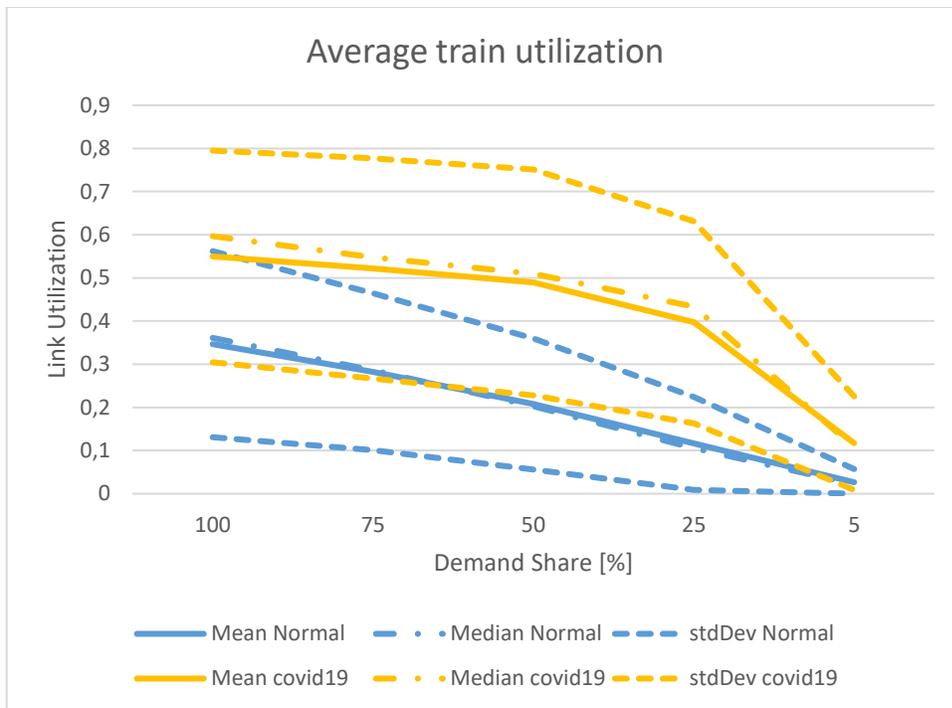

*Figure 6. Average train utilization in normal and covid19 conditions*

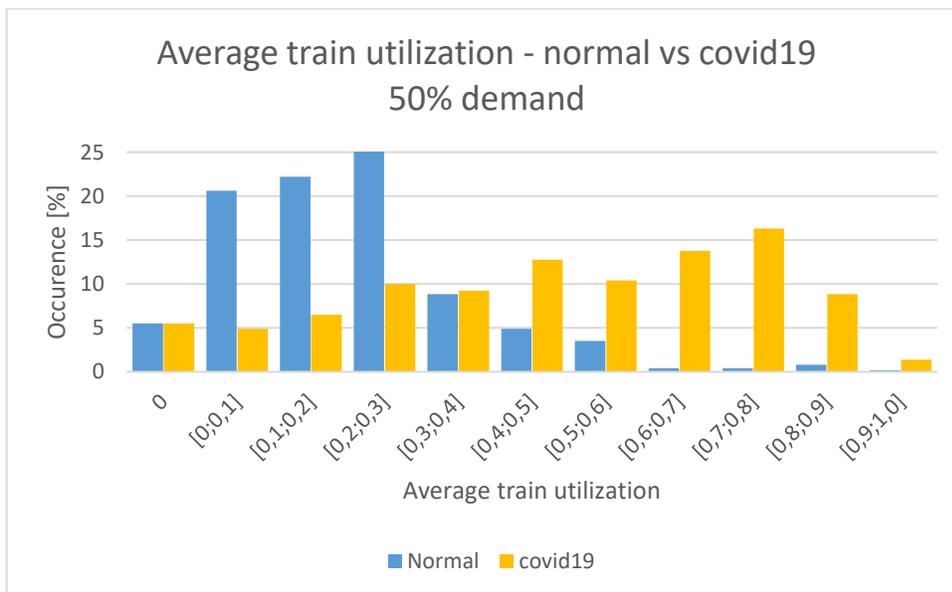

*Figure 7. Average train utilization in normal and covid19 for 50% demand*

Finally, **it is striking that for passengers demands of 25% of the original and higher, the system is expected to become more utilized based on its considered capacity than in normal conditions before covid19 as shown in Figure 8. Also, for 50% demand share every other train is maximally occupied at some point during its journey.** For normal conditions (resulting in average (light blue) and max values (dark blue)), we can state that a well-designed system performs with approximately equal mean and median values (this holds for both max and average train utilization). For example, with 100% of demand, max train utilization reaches 0.53 (mean) and 0.51 (median). Instead, for covid19, such balanced mean vs median ratio holds only for the demand shares of 5 and 25% while for bigger demand shares, the max utilization grows at a much higher rate. This represents that the system capacity becomes saturated pretty quickly. Therefore, it can be concluded **that the system under covid19 can perform "well" on the level of normal peak hours only up to 25% of normal demand.** In addition, looking at the ratio of trains with occupancy bigger than 0.90 at any link, these numbers are even more remarkable (Figure 9). In particular, **for covid19, about 67% of all trains would be least 90% full (max link util of 0.9) at some point of its journey**. While for normal, these numbers reach 0.20 only. Therefore, any demand of 25% and higher would mean experiencing higher train use and more crowding than in pre-covid19 conditions. This can clearly have significant implications on transport planning and management.

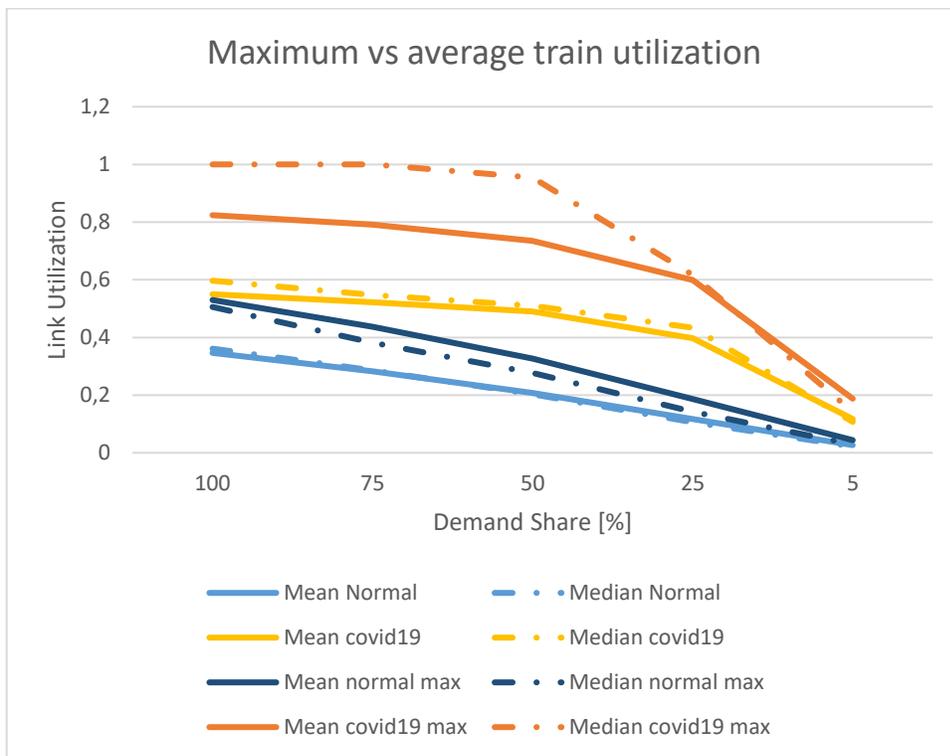

*Figure 8. Average vs max train utilization for normal and covid19 conditions*

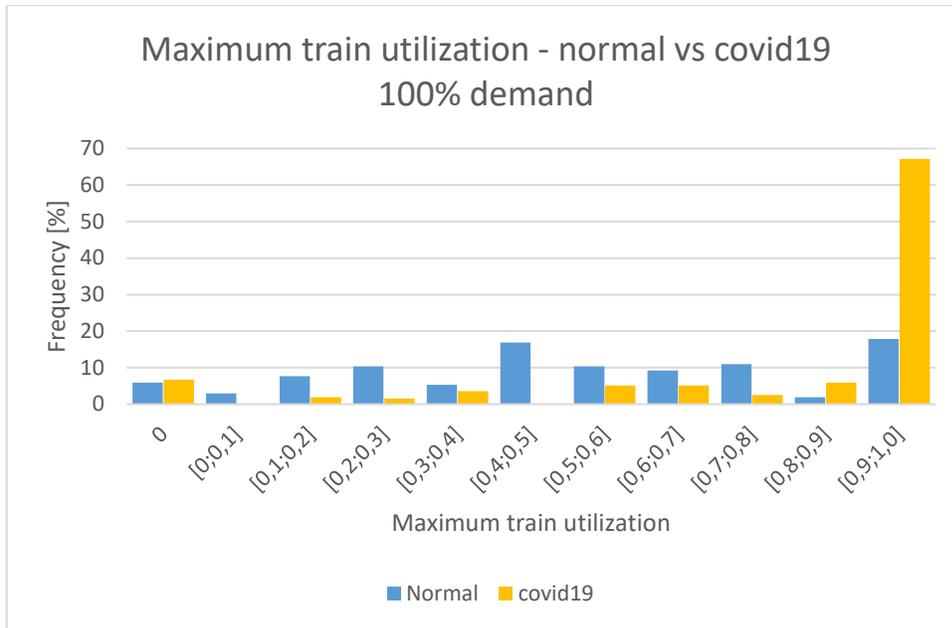

*Figure 9. Maximum train utilization in normal and covid19 for 100% demand*

Overall, for normal conditions, we observed linear changes in all 3 statistical indicators, mean, median and standard deviation, with changing the transport demand. For covid19 scenarios, we see a flattening behavior when increasing the demand (bigger than 25%), which suggests that trains (and) links are becoming saturated, i.e. reaching their offered/scheduled capacity. However, it may be expected that some links/trains do not reach their full capacity even when further increasing the demand in the current operations of railway system. There are two main reasons: such trains/links have only limited actual demand in its proximity, and they do not provide a satisfactory routing alternative to other passengers originating/ending further away. Finally, the **railway networks due to covid19 are likely to be more vulnerable to (multiple) disruptions**. In general, vulnerability of a network is mostly relieved by multiple good routing alternatives for passengers via less utilized links in case of disruptions, and an overall highly utilized network is likely to be lacking such alternatives. Thus, resulting from its high utilization and comparing to the normal case, the railway network during covid19 is expected to be significantly more vulnerable.

# 6 Conclusions and future research

We explored impacts of covid19 on transport capacity in the Dutch railway network using mathematical modelling. The results show that it can be expected that the current timetables during peak hours would not be sufficient with the restricted capacity due to physical distancing. In fact, the railway system could accommodate about 50% of the pre-covid19 demand. What is more, even with 25% of the normal demand not all passengers would be able to travel. And, capacity becomes significantly constraining as soon as the transport demand recovers towards 50% which leads to over 30% of not transported passengers and pretty full trains overall. Therefore, infrastructure managers and train operators shall be well prepared for these new challenges that can be expected in the very near future.

The presented model can be used for multiple purposes. First, it shall be useful to determine trains and areas that need more capacity to introduce extra train services on one side, and also the trains that are underutilized and might offer potential reallocation to lines where they are needed on the other side. Second, it can be used to evaluate a more detailed train usage (e.g. some passengers travel in pairs/groups, which may improve train capacity) and changes in network performance after applying various response measures. Third, it would be beneficial to understand expected exposure to and spreading of viruses in such railway and PT networks using epidemiological models, which would eventually help to understand what a "reasonable" train utilization is.

It becomes our reality that current railway system does not have sufficient transport capacity to satisfy the passenger demand under physical distancing due to covid19. Also, even though certain percentage of population may switch to home instead of office working, it may not be expected to have a significant impact on transport demand (anytime soon). We recognise several research directions that could support improving capacity usage and transport performance under these new conditions. First, new passenger-centered timetable (re-)design models that provide more transport capacity at most important links shall be addressed. Second, research on travel incentives to passengers and freight operators will be required. This may lead towards better capacity use of lightly utilized links/trains and overall to increased ridership. Third, researchers shall consider new concepts of travel reservation for rail transport as well as public transport in general. This would allow spreading the peak demand over a longer time. Fourth, this would lead to new traffic management requirements to manage both trains and passenger flows in railway and PT networks. Finally, new technological concepts shall be considered that allow more capacity on track such as virtual coupling, and alternative vehicle and station design could be considered.